\documentclass[11pt,twoside,reqno]{amsart}

\usepackage{amssymb}

\def\const{\text{\rm const}}

\def\sm{\setminus}
\def\ti{\tilde}

\def\dist{\text{\rm dist}}
\def\supp{\text{\rm supp}\,}
\def\to{\rightarrow}

\def\bs{\bigskip}

\def\ms{\medskip}
\def\no{\noindent}

\def\R{{\mathbb R}}

\def\Z{{\mathbb{Z}}}
\def\C{{\mathbb{C}}}
\def\N{{\mathbb{N}}}

\def\GG{{{{\mathbf{G}}}}}
\def\EE{{\mathcal E}}

\def\e{\varepsilon}
\def\d{\delta}
\def\L{\Lambda}
\def\l{\lambda}

\def\lan{\lambda_n}

\newcommand{\charf}{\raisebox{\depth}{\(\chi\)}}

\theoremstyle{plain}
\newtheorem{lem}{Lemma}
\newtheorem{thm}{Theorem}

\newtheorem{rem}{Remark}

\numberwithin{equation}{section}

\author{A.~Poltoratski}
\address{Texas A\&M University
\\ Department of Mathematics\\
College Station, TX 77843, USA }
\email{alexeip@math.tamu.edu}
\thanks{The author is supported by
N.S.F. Grant No. 0800300}

\title{A problem  on completeness of exponentials}

\begin{document}

\begin{abstract} Let $\mu$ be a finite positive measure on the real line.
For  $a>0$ denote by $\EE_a$ the family of exponential functions
$$\EE_a=\{e^{ist}| \ s\in[0,a]\}.$$
The exponential type of $\mu$ is the infimum of  all numbers $a$ such that the finite linear combinations of the exponentials from $\EE_a$ are dense in $L^2(\mu)$. If the set of such $a$ is empty, the exponential type of $\mu$ is defined as infinity. The well-known type problem asks to find the exponential type of $\mu$ in terms of $\mu$.

\ms\no In this note we present a solution to the type problem and discuss its relations with known results.
\end{abstract}

\maketitle

\ms
\section{\bf Introduction}

\ms\subsection{Completeness of exponentials.} Let $\mu$ be a finite positive Borel measure on $\R$. Let us consider
the family $\EE_\L$ of exponential functions $\exp(i\l t)$ on $\R$ whose frequencies $\l$ belong to a certain
set $\L\subset \C$:
$$\EE_\L=\{\exp(i\l t)|\ \l\in \L\}.$$
 One of the classical problems of Harmonic analysis is to find conditions on $\mu$ and $\L$ that ensure
 completeness, i.e. density of  finite linear combinations, of
  functions from $\EE_\L$  in $L^2(\mu)$.

 \ms\no Versions of this problem were considered by many prominent analysts.
 The case when $\L$ is a sequence and $\mu$ is Lebesgue measure on an interval was solved by Beurling and Malliavin
 in the early sixties \cite{BM1, BM2}. The so-called Beurling-Malliavin theory, created to treat that problem, is considered to be
 one of the deepest parts of the 20th century Harmonic Analysis.

 \ms\no Other cases of the problem and its multiple reformulations were studied by Wiener, Levinson, Kolmogorov, Krein and many others.
 Such an extensive interest is largely due to the fact that it is naturally related to other fields
 of classical analysis, such as stationary Gaussian processes and prediction theory, spectral problems for differential operators,
 approximation theory, signal processing, etc. Despite  considerable efforts by the analytic community
 many important cases of the problem remain open.

\ms
 \subsection{The type problem} Perhaps the most studied among such open cases is the so-called type problem.
 Consider a family $\EE_a=\EE_{[0,a]}$ of exponential functions whose frequencies belong to the interval
 from 0 to $a$.
 If $\mu$ is a finite positive measure on $\R$ we denote by $\GG^2_\mu$ its exponential type that is defined as
 \begin{equation} \GG^2_\mu=\inf\{\ a>0\ |\ \EE_a  \textrm{ is complete in }  L^2(\mu)\ \}\label{type}
 \end{equation}
 if the set of such $a$ is non-empty and as infinity otherwise. The type problem asks to calculate $\GG^2_\mu$ in terms
 of $\mu$.

 \ms\no This question first appears in the work of Wiener, Kolmogorov and Krein  in the context of stationary Gaussian processes (see
 \cite{Krein1, Krein2} or the book by Dym and McKean \cite{DM}). If $\mu$ is a spectral measure of a stationary Gaussian process, completeness of $\EE_a$ in $L^2(\mu)$
 is equivalent to the property that  the process at any time is determined by the data for the time
 period from 0 to $a$. Hence the type of the measure is the minimal length of the period of observation necessary to predict the rest of the process. Since any even measure is a spectral measure of a stationary Gaussian process, and vice versa, this reformulation is practically equivalent.

 \ms\no The type problem can also be restated in terms of the Bernstein weighted approximation, see for instance the book by Koosis
 \cite{KoosisLog}. Important connections with spectral theory of second order differential operators were studied by Gelfand and Levitan \cite{GL} and Krein \cite{Krein2, Krein3}.

\ms\no Closely related to spectral problems for differential operators is Krein -- de Branges' theory of Hilbert spaces of entire functions, see \cite{dBr}.
One of the deep results of the theory says that for any positive finite (or more generally Poisson-finite) measure $\mu$ on $\R$ there
is a unique nested regular chain of de Branges' spaces of entire functions isometrically embedded in $L^2(\mu)$. An important characteristic of such a chain
is the supremum $S_\mu$ of the exponential type taken over all entire functions contained in the embedded spaces. For instance, if such a chain corresponds to a
regular Schr\"odinger operator on an interval, i.e. if  $\mu$ is the spectral measure of such an operator, then $S_\mu$ is equal to the length of the interval and all spaces of the chain can be parametrized by their exponential type. It is well-known, and not difficult to show,
that the problem of finding the value of $S_\mu$ is equivalent to the type problem, i.e. $S_\mu=G^2_\mu$.

 \ms\no For more on the history and connections of the type problem see, for instance,  a note by Dym \cite{Dym} or a recent paper by Borichev and Sodin \cite{BS}.

\ms
 \subsection{General case $p\neq 2$}
 The family $\EE_a$ is incomplete in $L^2(\mu)$ if and only if there exists a function $f\in L^2(\mu)$ orthogonal to all
 elements of $\EE_a$. Expanding  to other $1\leqslant p\leqslant \infty$ we define

 \begin{equation} \GG^p_\mu=\sup\{\ a\ |\  \exists \ f\in L^p(\mu),\int f(x)e^{i\l x}d\mu(x)=0,  \forall \ \l\in[0,a]\ \}.
 \label{typep}\end{equation}

 \ms\no We put $\GG^p_\mu=0$ if the set in \eqref{typep} is empty. By duality, for $1<p\leqslant\infty$, $\GG^p_\mu$ can still be defined as the infimum of $a$ such that $\EE_a$ is complete in $L^q(\mu),\ \frac 1p +\frac 1q =1$.
Cases $p\neq 2$ were considered in several papers, see for instance articles by Koosis \cite{Koosis2} or Levin \cite{Levin2} for the case $p=\infty$ or \cite{GAP} for $p=1$.

\ms\no Since $\mu$ is a finite measure we have
\begin{equation}\GG^p_\mu\leqslant\GG^q_\mu\textrm{ for }p\geqslant q.\label{pq}\end{equation} Apart from this obvious observation, the problems of
finding $\GG^p_\mu$ for different $p$ were generally considered non-equivalent. One of the consequences of
theorem \ref{main}, section \ref{secMain}, is that, in some sense, there are only two significantly different cases, $p=1$ (the gap problem)
and $1<p\leqslant\infty$ (the general type problem).

\ms\no In this paper we restrict our attention to the class of finite measures. The formal reason for that is the fact that $\mu$ has
to be finite for exponentials to belong to $L^2(\mu)$. This obstacle can be easily overcome if instead of $\EE_a$ one considers $E_a$,
the set of Fourier transforms of smooth functions supported on $[0,a]$. All elements of $E_a$ decay fast at infinity and
one one can ask about the density of $E_a$ in $L^p(\mu)$ for wider sets of $\mu$, see for instance \cite{BS}. One of such traditional sets is the class of Poisson-finite measures satisfying
$$\int\frac {d|\mu|(x)}{1+x^2}<\infty.$$
However, due to the reasons similar to lemma \ref{polynomial} below
(note that if $\mu$ is Poisson-finite then $\mu/(1+x^2)$ is finite and
vice versa), considering such a wider set of measures will not change the problem
and all of the statements will remain the same or analogous.

\ms
\subsection{The gap problem}\label{introGAP}
One of the important cases is the so-called gap problem, $p=1$. Here one can reformulate the question
as follows.

\ms\no Let $X$ be a closed subset of the real line. Denote
$$\GG_X=\sup\{\ a\ |\ \exists\   \mu\neq 0, \ \supp\mu\subset X,   \hat\mu=0   \textrm{ on }[0,a]\ \}.$$
Here and in the rest of the paper $\hat\mu$ denotes the (inverse) Fourier transform of a finite measure $\mu$ on $\R$:
$$\hat\mu(z)=\int_\R e^{izt}d\mu(t).$$
As was shown in \cite{GAP}, for any finite measure $\mu$ on $\R$, $\GG^1_\mu$, as defined in the previous section, depends only on its support:
$$\GG^1_\mu=\GG_X,\ X=\supp\mu.$$
This property separates the gap problem from all the cases $p>1$.

\ms\no For a long time both the gap problem and the type problem were considered
 by experts to be "transcendental," i.e. not having a closed form solution.
Following an approach developed in \cite{MIF1} and \cite{MIF2},
a solution to the gap problem was recently suggested in \cite{GAP}, see section \ref{secGAP}. Some of definitions and results from \cite{GAP} are used in the present paper.

\ms
\subsection{Known examples} We say that a function $f$ on $\R$ is Poisson-summable if it is summable with respect to the Poisson measure $\Pi$,
$$d\Pi=dx/(1+x^2).$$
We say that a sequence of real numbers $A=\{a_n\}$ is discrete if
it does not have finite accumulation points.
We always assume that a discrete sequence is enumerated
in the natural increasing order:  $a_n\geqslant a_{n-1}$. Since the sequences considered here
have $\pm\infty$ as their density points, the indices run over $\Z$.
In most of our statements and definitions, the sequences do not have multiple points.
 We call a discrete sequence $\{a_n\}\subset \R$ separated if $|a_n-a_k|>c$ for some $c>0$ and any $n\neq k$.

\ms\no A classical result by Krein \cite{Krein1} says that if $d\mu=w(x)dx$ and $\log w$ is Poisson-summable then
$\GG^p_\mu=\infty$ for all $p, \ 1\leqslant p\leqslant\infty$. A partial inverse, proved by Levinson and McKean, holds for even monotone $w$,
see section \ref{secKLM}.

\ms\no A theorem by Duffin and Schaeffer \cite{DS} implies that if $\mu$ is a measure such that
for any $x\in \R$
$$\mu([x-L,x+L])>d$$
for some $L,d>0$ then $\GG^2_\mu\geqslant 2\pi/L$, see section \ref{secDS}.

\ms\no For discrete measures, in the case $\supp\mu=\Z$,  a deep result by Koosis shows an
analogue of Krein's result: if $\mu=\sum w(n)\delta_n$, where
$$\sum\frac{\log w(n)}{1+n^2}>-\infty,$$
then $\GG^p_\mu=2\pi$ for all $p, \ 1\leqslant p\leqslant\infty$ \cite{Koosis2}. Not much was known about supports other than $\Z$ besides a recent result from \cite{Polya}, which implies that if
$$\mu=\sum \frac{\delta_{a_n}}{1+a_n^2}$$
for a separated sequence $A=\{a_n\}\subset\R$ then $\GG^p_\mu=2\pi D_*(A)$, where $D_*$ is the interior Beurling-Malliavin density of $A$, see section \ref{secGAP} for the definition. We generalize these results in section \ref{secDisc}.

\ms\no In addition to these few examples, classical theorems by Levinson-McKean, Beurling and de Branges
show that if a measure has long gaps in its support or decays too fast, then $\GG^p_\mu=0$, see section \ref{classical}.
Examples of measures of positive type  can  be constructed using the results by Benedicks \cite{Benedicks}, see section \ref{secBen}.
The most significant recent development, that allows one to modify existing examples, is the result by Borichev and Sodin
\cite{BS},
which says that "exponentially small" changes in weight or support do not change the type of a measure, see section
\ref{secBS}.

\ms
 \subsection{Approach and goals of the paper}
 The problems discussed above belong to the area  often called the Uncertainty Principle in Harmonic Analysis \cite{HJ}.
 A new approach developed by N. Makarov and the author in \cite{MIF1, MIF2} allows one to study
 this area with modern tools of analytic function theory and singular integrals. Together with traditional methods, such as de Branges'
 theory of Hilbert spaces of entire functions or the Beurling-Malliavin theorems, these techniques have produced some new
  ideas and developments. Among them is
 an extension of the Beurling-Malliavin theory \cite{MIF2}, a solution to the P\'olya-Levinson problem on
 sampling sets for entire functions of zero type \cite{Polya} and a solution to the gap problem \cite{GAP}.
 In the present paper we continue to apply the same approach.

 \ms\no We focus on the type problem, the problem
 of finding $\GG^2_\mu$ in terms of $\mu$.
 Our main results are theorem \ref{main} and its corollaries contained in section \ref{mainresults}.
 In most of our statements, treating  $p>1,\ p\neq 2$ did not require
 any additional efforts, and hence they were formulated for general $p>1$.
 The case $p=1$, studied in \cite{GAP}, provided us with some useful definitions
 and statements, see section  \ref{secGAP}.

\ms \textbf{Acknowledgements.} I am grateful to Nikolai Makarov whose deep mathematical insight and intuition led to the development of the methods used in this paper. I would also like to thank Misha Sodin for getting me interested in the gap and type problems and for numerous invaluable discussions.

\ms
 \subsection{Contents}
 The paper is organized as follows:

\begin{itemize}
\item Section \ref{prelim}  contains preliminary material, including the basics of the so-called Clark theory, definitions
of Beurling-Malliavin densities and a short discussion of the gap problem.

\item In section \ref{mainresults} we state the main results of the paper.

\item Section \ref{classical} discusses connections of our results with classical theorems by  Beurling,
de Branges, Duffin and Schaeffer, Krein, Levinson and McKean as well as more recent results by Benedicks, Borichev and Sodin.

\item Section \ref{lemmas} contains  several lemmas needed for the main proofs.

\item In section \ref{mainproofs} we give the proofs of the main results.
\end{itemize}

 \ms\section{Preliminaries}\label{prelim}


\ms
\subsection{Clark theory}\label{Clark}

\ms\no By $H^2$ we denote the Hardy space in the upper half-plane
$\C_+$. We say that an inner function
$\theta(z)$ in $\C_+$ is  meromorphic if it allows a
meromorphic extension to the whole complex plane. The meromorphic
extension to the lower half-plane $\C_{-}$ is given by
$$\theta(z)=\frac{1}{\theta^{\#}(z)} $$
where $\theta^{\#}(z)=\bar\theta(\bar z)$.

\ms\no Each inner function $\theta(z)$ determines a model subspace
$$K_\theta=H^2\ominus \theta H^2$$ of the Hardy space
$H^2(\C_+)$. These subspaces play an important role in complex and
harmonic analysis, as well as in operator theory,
see~\cite{Ni2}.

\ms\no For each inner function $\theta(z)$ one can consider a positive harmonic
function
$$\Re \frac{1+\theta(z)}{1-\theta(z)}$$
 and, by the Herglotz
representation, a positive measure $\mu$ such that
\begin{equation} \label{for1} \Re
\frac{1+\theta(z)}{1-\theta(z)}=py+\frac{1}{\pi}\int{\frac{yd\mu
(t)}{(x-t)^2+y^2}}, \hspace{1cm} z=x+iy,\end{equation} for some $p
\geqslant 0$. The number $p$ can be viewed as a point mass at infinity.
The measure $\mu$ is Poisson-finite, singular and supported on the set where non-tangential limits
of $\theta$ are equal to $1$.
The measure
$\mu +p\delta_\infty$ on $\hat\R$
is called the Clark measure for $\theta(z)$.

\ms\no Following  standard notations, we will sometimes denote the Clark measure defined in \eqref{for1} by $\mu_1$.
More generally, if $\alpha\in \C, |\alpha|=1$ then $\mu_\alpha$ is the measure defined by \eqref{for1}
with $\theta$ replaced by $\bar\alpha\theta$.

\ms\no Conversely, for
every positive singular Poisson-finite measure $\mu$  and a number $p \geqslant
0$, there exists an inner function $\theta(z)$ satisfying~\eqref{for1}.

\ms\no Every function $f \in K_\theta$ can be represented by the formula
\begin{equation} \label{for2} f(z)=\frac{p}{2\pi
i}(1-\theta(z))\int{f(t)\overline{(1-\theta(t))}dt}+\frac{1-\theta(z)}{2\pi
i}\int{\frac{f(t)}{t-z} d\mu (t)}.
\end{equation}
If the Clark measure does not have a point mass at infinity,
the formula is simplified to
$$f(z)=\frac1{2\pi i}(1-\theta(z))Kf\mu$$
where $Kf\mu$ stands for the Cauchy integral
$$Kf\mu(z)=\int\frac{f(t)}{t-z} d\mu (t).$$
This gives an isometry of
$L^2(\mu)$ onto $K_\theta$.
Similar formulas can be written for any $\mu_\alpha$ corresponding to $\theta$.

\ms\no In the case of meromorphic $\theta(z)$,
every function $f \in K_\theta$ also has a meromorphic extension in
$\C$, and it is given by the formula~\eqref{for2}. The
corresponding Clark measure is discrete with atoms at the points
of $\{\theta=1\}$ given by
$$\mu(\{x\})=\frac{2\pi
}{|\theta'(x)|}.$$
If $\L\subset \R$ is a given discrete sequence, one can easily construct a meromorphic
inner function $\theta$ satisfying $\{\theta=1\}=\L$ by  considering a positive
Poisson-finite measure concentrated on $\L$ and then choosing $\theta$ to satisfy
\eqref{for1}. One can prescribe the derivatives of $\theta$ at $\L$ with a proper choice
of pointmasses.

\ms\no For more details on Clark measures and further references the reader may consult \cite{PS}.

 \ms\subsection{Interior and exterior densities}
A sequence of disjoint intervals $\{I_n\}$ on the real line is called \textit{long} (in the sense of Beurling and Malliavin)
if
\begin{equation}\sum_n\frac{|I_n|^2}{1+\dist^2(0,I_n)}=\infty,\label{long}
\end{equation}
where $|I_n|$ stands for the length of $I_n$.
If the sum is finite, we call $\{I_n\}$ \textit{short}.

\ms\no One of the obvious properties of short sequences  is that
$|I_n|=o(\dist(0,I_n))$ as $n\to\infty$.
In particular, $\dist(0,I_n)$ can be replaced with any $x_n\in I_n$ in \eqref{long}.

\ms\no Following \cite{BM2} we say that a discrete sequence $\L\subset
\R$ is $a$-\textit{regular} if for every $\epsilon>0$ any sequence of
disjoint intervals $\{I_n \}$ that satisfies
$$\left|\frac{\#(\L\cap I_n)}{|I_n|}-a\right|\geqslant
\epsilon$$
for all $n$, is short.

\ms\no A slightly different $a$-regularity can be defined in the following way,
that is more convenient in some settings.
For a discrete sequence $\L\subset \R$ we denote by $n_\L (x)$ its
 counting function, i.e. the step function  on $\R$, that is constant between any two points   of $\L$, jumps up by $1$ at each point of $\L$ and is equal to $0$ at $0$.
We say that $\L$ is \textit{strongly} $a$-\textit{regular}
if
$$\int\frac{|n_{\L}(x)-ax|}{1+x^2}<\infty.$$

\ms\no Conditions like this can be found in many related results, see for instance
\cite{dBr} or \cite{KoosisLog}. Even though $a$-regularity is not equivalent to
strong $a$-regularity, in the following definitions of densities changing
"$a$-regular" to "strongly $a$-regular"  will lead to equivalent definitions.

\ms\no The interior BM (Beurling-Malliavin)
density of a sequence $\L$ is defined as
\begin{equation}\label{id}D_*(\L):=\sup \{\ a\ |\ \exists \ \textrm{$a$-regular subsequence}\ \ \L'\subset \L \ \}.\end{equation}
If the set is empty we put $D_*(\L)=0$.
Similarly, the exterior BM density is defined as
\begin{equation}\label{ed}D^*(\L):=\inf \{\ a\ |\ \exists \ \textrm{$a$-regular supsequence}\ \ \L'\supset \L\ \}.\end{equation}
\no If no such sequence exists, $D^*(\L)=\infty$.

\ms\no It is interesting to observe that after the two densities were  simultaneously introduced
over fifty years ago, the exterior density immediately became one of the staples of harmonic
analysis and spectral theory, mostly due to its appearance in the celebrated Beurling-Malliavin theorem,
see \cite{BM2}, \cite{HJ} or \cite{KoosisLog}.
Meanwhile, the interior density remained largely forgotten until its recent comeback in
\cite{Polya} and \cite{GAP}. It will continue to play an important role in our discussions below.

\ms\subsection{The gap problem and $d$-uniform sequences}\label{secGAP}
Let $\Lambda=\{\lambda_1,...,\lambda_n\}$ be a finite set of points on $\R$. Define

\begin{equation}E(\Lambda)=\sum_{\lambda_k,\lambda_j\in\L} \log|\lambda_k-\lambda_l|.\label{electrons}\end{equation}
According to the 2D Coulomb law, the quantity $E(\L)$ can be interpreted as potential energy of the system of "flat electrons" placed at $\L$,
see \cite{GAP}. That observation motivates the term we use for the condition \eqref{energy} below.
\bs

\ms\no The following example is included to illustrate our next definition.

\bs\no \textbf{Key example:}

\ms\no \textit{Let $I\subset \R$ be an interval and
let $\L=d^{-1}\Z\cap I$ for some $d>0$.
Then
$$\Delta=\#\L=
d|I|+O(1)$$
and
\begin{equation}E=E(\Lambda)=\sum_{1\leqslant m\leqslant \Delta}\log\left[d^{-\Delta+1}(m-1)!(\Delta-m)!\right]=
\Delta^2\log |I| +  O(|I|^2)\label{eqkey}\end{equation}
as follows from Stirling's formula. Here  the notation $O(\cdot)$ corresponds to the direction $|I|\to\infty$.
}

\begin{rem}
The uniform distribution of points on the interval does not maximize
the energy $E(\Lambda)$ but comes within $O(|I|^2)$ from the maximum, which is negligible
 for our purposes, see the main definition and its discussion below. It is interesting to observe that the
maximal energy for $k$ points is achieved when the points are placed at the endpoints of $I$ and the zeros of
the Jacobi $(1,1)$-polynomial of degree $k-2$, see for example \cite{Kerov}.
\end{rem}

\ms\no Let
$$...<a_{-2}<a_{-1}<a_0=0<a_1<a_2<...$$
 be a discrete sequence of real points. We say that the intervals $I_n=(a_n,a_{n+1}]$ form a short partition
of $\R$ if $|I_n|\to\infty$  as $ n\to \pm\infty$ and the sequence $\{I_n\}$ is short.

\bs\no \textbf{Main Definition:}

\ms\no Let $\L=\{\lan\}$ be a discrete sequence of real points. We say that $\L$ is
$d$-uniform if
there exists a short partition $\{I_n\}$ such that
\begin{equation} \Delta_n= d|I_n|+o(|I_n|)\ \ \   \text{for all} \ \ \ n\ \ \textrm{(density condition)}\label{density}\end{equation}
as $n\to\pm\infty$ and
\begin{equation} \sum_n \frac{\Delta_n^2\log|I_n|-E_n}{1+\dist^2(0,I_n)}<\infty\ \ \textrm{(energy condition)}\label{energy}\end{equation}
where $\Delta_n$ and $E_n$ are defined as
$$\Delta_n=\#(\L\cap I_n)\ \ \text{ and }\ \ E_n=E(\Lambda\cap I_n)=\sum_{\lambda_k,\lambda_l\in I_n,\ \lambda_k\neq\lambda_l}\log|\lambda_k-\lambda_l|.$$

\ms\begin{rem}
Note that the series in the energy condition is positive:  every term in the sum defining $E_n$ is at most $\log|I_n|$ and there
are  no more than $\Delta_n^2$ terms.

\ms\no As follows from the example above, the first term in the numerator of \eqref{energy} is approximately equal to
the energy of $\Delta_n$ electrons spread uniformly over $I_n$. The second term is the energy of electrons placed at
$\L\cap I_n$. Thus the  energy condition
is a requirement that the placement of the points of $\Lambda$ is close to uniform, in the sense that the work needed to
spread the points of $\L$ uniformly on each interval is
summable with respect to the Poisson weight. For a more detailed discussion of this definition see \cite{GAP}
\end{rem}

\ms\no In \cite{GAP}, $d$-uniform sequences were used to solve the gap problem mentioned in the introduction.
Recall that with any closed $X\subset\R$ one can associate its (spectral) gap characteristic $\GG_X$ defined as
in section \ref{introGAP}.
The main result of \cite{GAP} is the following statement:

\begin{thm}\label{mainGAP}\cite{GAP} Let $X$ be a closed set on $\R$. Then
$$\GG_X=\sup\{\ d\ |\  X \textrm{ contains a }
d-\textrm{uniform sequence }\}.$$
\end{thm}

\no Recall that, as was proved in \cite{GAP}, $\GG_X=\GG^1_\mu$ for any $\mu$ such that $\supp\mu=X$.
The following simple observations will also be useful to us in the future:

\ms
\begin{rem}\label{rem1}

$$$$

\begin{itemize}

\item If $\L$ is a $d$-uniform sequence then $D_*(\L)=d$, as follows easily from  the density condition \eqref{density}.

\item Among other things, the energy condition ensures that the points of $\L$ are not too close to each other. In particular, if
$\L$ is $d$-uniform for some $d>0$ and
$\L'=\{\lambda_{n_k}\}$ is a subsequence such that for all $k$,
$$\lambda_{n_k+1}-\lambda_{n_k}\leqslant e^{-c|\lambda_{n_k}|}$$
for some $c>0$, then $D_*(\L')=0$.

\item An exponentially small perturbation of a $d$-uniform sequence contains a $d$-uniform subsequence. More precisely,
if $c>0$ and $\L$ is a $d$-uniform sequence then any sequence $A=\{\alpha_n\}$ such that
$|\lan-\alpha_n|\leqslant e^{-c|\lan|}$ contains a $d$-uniform subsequence $A'$ consisting of
 all $\alpha_{n_k}$ such that
$$\lambda_{n_k+1}-\lambda_{n_k}\geqslant e^{-(c-\e)|\lambda_{n_k}|}.$$

\item As discussed in \cite{GAP}, the energy condition always holds for separated sequences. If $\L$ is separated
then it is $d$-uniform if and only if $D_*(\L)=d$.

\end{itemize}
\end{rem}

\ms\subsection{Polynomial decay}\label{poly}
In this section we prove a version of the well-known property that adding or removing polynomial decay cannot change the type of
a measure.

\begin{lem}\label{polynomial}
Let $\mu$ be a finite positive measure on $\R$ and let $\alpha>0$. Consider the measure $\nu$ satisfying
$$d\nu(x)=\frac{d\mu(x)}{1+|x|^\alpha}.$$
Then for any $1\leqslant p\leqslant \infty$
$$\GG^p_\mu=\GG^p_\nu.$$

\end{lem}

\begin{proof}
Since $d\nu/d\mu\leqslant 1$, one only needs to show that $\GG^p_\mu\leqslant \GG^p_\nu.$
Suppose that $f\in L^p(\mu)$ is such that $\bar f\mu$ annihilates all $e^{iaz}, a\in (0,d)$.
This is equivalent to the property that the Cauchy
integral $Kf\mu$ is divisible by $e^{idz}$ in $\C_+$, i.e. it decays like $e^{idz}$ along the positive imaginary axis $i\R_+$, see for instance lemma 2 in \cite{Polya}.

\ms\no Let $N\geqslant\alpha$ be an integer. It is enough to prove the statement for $N=1$: the general case will follow
by induction.

\ms\no First let us assume that $Kf\mu$ has at least one zero $a$ in $\C\setminus\R$. It is well-known, and not difficult
to verify, that then the measure $\frac f{x-a}\mu$ satisfies
$$K\left(\frac f{x-a}\mu\right)=\frac {Kf\mu}{z-a}.$$
Hence the Cauchy integral in the left-hand side  still decays like $e^{idz}$ along $i\R_+$ and therefore the measure still annihilates $e^{iaz}, a\in (0,d)$. It is left to notice that
$$f(x)\frac{1+|x|^\alpha}{x-a}\in L^p(\nu).$$

\ms\no If $Kf\mu$ does not have any zeros outside of $\R$, note that the Cauchy integral of the measure $\eta=e^{-i\e x}f\mu$ satisfies
$$K\eta=e^{-i\e z}Kf\mu$$
(see for instance theorems 3.3 and 3.4 in \cite{Clark}) and therefore
$$K(f\mu-c\eta)=K(1-ce^{-i\e x})f\mu$$
has infinitely many zeros in  $\C\setminus\R$ for any $c, |c|\neq 1$, while still
decaying like $e^{i(d-\e)z}$ along $i\R_+$.
\end{proof}

\ms\no Let $\mu$ be a finite measure on a separated sequence $X$, with point masses decaying polynomially.
Lemma \ref{polynomial} together with  elementary estimates imply that then $\GG^2_\mu=\GG^1_\mu$.
Hence in this case theorem \ref{main} below becomes the statement from \cite{Polya} mentioned above:
$$\GG^2_\mu=2\pi D_*(X).$$

\ms
\section{Main Results}\label{mainresults}

\ms
\subsection{Main Theorem}\label{secMain}
Let $\tau$ be a finite positive  measure on the real line. We say that
a  function $W\geqslant 1$ on $\R$ is a $\tau$-\textit{weight} if $W$ is lower semi-continuous,
tends to $\infty$ at $\pm\infty$  and $W\in L^1(\tau)$.

\begin{thm}\label{main}
Let $\mu$ be a finite positive measure on the line. Let $1<p\leqslant\infty$ and $a>0$ be  constants.

\ms\no Then $\GG^p_\mu\geqslant a$
if and only if for any  $\mu$-weight $W$ and any $0<d<a$ there exists
a $d$-uniform sequence $\L=\{\lan\}\subset \supp \mu$ such that
\begin{equation}\sum\frac{\log W(\lan)}{1+\lan^2}<\infty.\label{ur4}\end{equation}

\end{thm}

\ms\no We postpone the proof until section \ref{mainproofs}.

\ms\no One of the immediate corollaries of the above statement is that the $p$-type of a measure,
$\GG^p_\mu$ for $1<p\leqslant\infty$, does not depend on $p$, which may come as a surprise to some
of the experts.
Further corollaries of theorem \ref{main} and its connections with classical results are discussed
in the following sections.

\bs\subsection{Discrete case}\label{secDisc} The conditions of theorem \ref{main} are simplified for many specific classes of measures.
In particular, if the measure is discrete, or absolutely continuous with regular enough density, the weight $W$ may be eliminated from
the statement. Here we treat the discrete case that is important in spectral theory of differential operators and other adjacent areas.
Our results in this section may be viewed as extensions of the result by Koosis mentioned in the introduction.

\ms\no The following statement gives a simplified formula for the type of a measure supported on a discrete sequence, excluding pathological
cases when the counting function of the sequence grows exponentially.

\begin{thm}\label{mainDiscrete} Let $B=\{b_n\}$ be a discrete sequence of real points. Let  $$\mu=\sum w(n)\delta_{b_n}$$
be a finite positive measure supported on $B$.
Define
$$D = \sup\{\ d\ |\ \exists\textrm{ d-uniform }B'\subset B,\ \sum_{\lan\in B'}\frac {\log w(n)}{1+{n}^2}>-\infty\ \}.$$
Then for any $1<p\leqslant\infty$,
$$\GG^p_\mu\geqslant 2\pi D.$$
If the counting function of $B$ satisfies $\log (|n_B|+1)\in L^1_\Pi$ then
$$\GG^p_\mu=2\pi D.$$

\end{thm}

\ms\no The proof is given in section \ref{mainproofs}. The condition $\log (|n_B|+1)\in L^1_\Pi$ in the second part of the statement is
sharp. The corresponding examples can be easily constructed using theorem \ref{main} or the result by Borichev and Sodin \cite{BS},
see theorem \ref{BSthm} below.

\ms\no In the case when the sequence is separated, the condition can be simplified even further.
Note that for $p=1$, $\GG^1_\mu=D_*(\L)$ for any separated sequence $\L$ and any measure $\mu,\ \supp\mu=\L$,
by theorem \ref{mainGAP}. For $p>1$ we have

\begin{thm}\label{mainSeparated} Let $\L=\{\lambda_n\}$ be a separated sequence and let
$$\mu=\sum w(n)\delta_{\lambda_n}$$ be a finite positive measure supported on $\L$.
Define $$D=\sup D_*(\L'),$$
where the supremum is taken over all subsequences $\L'\subset \L$ satisfying
\begin{equation}\sum_{\lambda_n\in \L'}\frac{\log w(n)}{1+n^2}>-\infty.\label{ur10}
\end{equation}
Then $$\GG^p_\mu=2\pi D$$
for all $1<p\leqslant \infty$.
\end{thm}

\begin{proof}
Suppose that $G^p_\mu>2\pi D$ for some $D>0, p>1$.
Define the $\mu$-weight $W$ as $W(\lambda_n)=(\mu(\{\lambda_n\})(1+\lambda^2_n))^{-1}$.
Then by theorem \ref{mainDiscrete} there
exists a subsequence $\L'\subset \L$ such that $D_*(B')>D$ and \eqref{ur10} is
satisfied.

In the opposite direction the statement follows directly from theorem
\ref{mainDiscrete} and remark \ref{rem1}.
\end{proof}

\bs

\ms\subsection{A general sufficient condition}
As a corollary of theorem \ref{mainDiscrete} we obtain the following sufficient condition for general measures.
The condition seems to be reasonably sharp, as it is satisfied by all  examples of measures with positive type
existing in the literature.

\begin{thm}\label{SuffGen}
 Let $\mu$ be a finite positive measure on $\R$. Let
$A=\{a_n\}$ be a $d$-uniform sequence of real numbers such that
\begin{equation}\sum\frac{\log\mu((a_n-\e_n,a_n+\e_{n}))}{1+n^2}>-\infty,\label{log1}\end{equation}
 where
 $$\e_n=\frac 13\min\left((a_{n+1}-a_n),(a_n-a_{n-1})\right).$$

\ms\no Then
$\GG^\infty_\mu\geqslant 2\pi d.$
\end{thm}

\begin{proof} For each $\tau\in [0,1]$ let us define a discrete measure $\nu_\tau$ as follows.
The measure $\nu_\tau$ has exactly one pointmass of the size
$$\mu((a_n-\e_n,a_n+\e_{n}))$$
 in each interval
$$(a_n-\e_n,a_n+\e_{n})$$
at the point $x^\tau_n$ chosen as
$$x^\tau_n=\inf\{\ a\ |\ \mu((a_n-\e_n,a))\geqslant\tau \mu((a_n-\e_n,a_n+\e_{n}))\ \}.$$
Notice that $\{x^\tau_n\}$ is a $d$-uniform sequence. In view of \eqref{log1} and theorem \ref{mainDiscrete},
$\nu_\tau$
 satisfies
 $$\GG^\infty_{\nu_\tau}\geqslant 2\pi d.$$
  Then
  $$\nu=\int_0^1\nu_\tau d\tau$$
   satisfies $d\nu/d\mu\leqslant 1$ and therefore
$$\GG^\infty_\mu\geqslant\GG^\infty_\nu\geqslant 2\pi d.$$
\end{proof}

\section{Classical results and further corollaries}\label{classical}

\ms\no The goal of this section is to give examples of applications of theorem \ref{main}
and discuss its connections with classical  results on the type problem.
Due to this reason, we prefer to deduce each statement directly from the results of the last section,
rather than obtaining them from each other, even when the latter approach may slightly shorten the proof.

\ms\no In our estimates we write $a(n)\lesssim b(n)$ if $a(n)< C b(n)$ for some positive
constant $C$, not depending on $n$, and large enough $|n|$. Similarly, we write $a(n)\asymp b(n)$ if $c a(n)< b(n)<C a(n)$
for some $C\geqslant c>0$. Some formulas will have other parameters in place of $n$ or no parameters at all.

\ms
\subsection{Beurling's Gap Theorem}
\begin{thm}[Beurling \cite{Beurling-Stanford}]
If $\mu$ is a finite measure supported on a set with long gaps and
the Fourier transform of $\mu$ vanishes on an interval, then $\mu\equiv 0$.

\end{thm}

\begin{proof}
If $\supp\mu$ has long gaps than for every short partition of $\R$ infinitely many intervals
of the partition must be contained in the gaps of $\supp\mu$. Therefore $\supp\mu$ does
not contain a sequence satisfying the density condition \eqref{density}, i.e. it does not
contain a $d$-uniform sequence for any $d>0$.
\end{proof}

\ms\subsection{Levinson's Gap Theorem}

\begin{thm}[Levinson \cite{Levinson}]\label{Lev}
Let $\mu$ be a finite measure on $\R$ whose Fourier transform
vanishes on an interval. Denote
$$M(x)=|\mu|((x,\infty)).$$
If $\log M$ is not Poisson-summable on $\R_+$ then $\mu\equiv 0$.
\end{thm}

\begin{proof} Suppose that $\log M$ is not Poisson-summable on $\R_+$.
Without loss of generality, $M(0)=1$. Let $0=a_0\leqslant a_1 \leqslant a_2\leqslant...$
be the sequence of points such that
$$a_n=\inf\{\ a\ |\ M(a)\leqslant 3^{-n}\ \}.$$
 Define a $|\mu|$-weight $W$ as $2^n$ on each $(a_{n-1},a_n],\ a_{n-1}<a_n$.

\ms\no Since $\hat\mu$ vanishes on an interval, by theorem \ref{main} there exists a sequence $\L\subset\supp\mu$
 satisfying the density condition \eqref{density} with some $a>0$ on a short partition $I_n=(b_n,b_{n+1}]$, such that \eqref{ur4}
 holds.  WLOG $b_0=0$.
Notice that $\log W$ is an increasing step function on $\R_+$
satisfying $\log W\gtrsim -\log M$. Also, since $\{I_n\}$ is short,
$cb_{n+1}\leqslant  b_n$ for some $0<c<1$ and all $n>0$. Hence,
$$\sum_n\frac{\log W(\lan)}{1+\lan^2}\gtrsim\sum_{n=1}^\infty\frac{\log W(b_n)|I_n|}{1+b_n^2}$$$$
\gtrsim \sum_{n=1}^\infty\frac{\log W(cb_{n+1})|I_n|}{1+b_n^2}\gtrsim \int_0^\infty \frac{-\log M(cx) dx}{1+x^2}=\infty.$$

\end{proof}

\ms\no Levinson's result above was later improved by Beurling \cite{Beurling-Stanford} who showed that instead of
vanishing on an interval $\hat\mu$ may vanish on a set of positive Lebesgue measure
with the same conclusion. Note that an analogous improvement cannot be made in Beurling's own
gap theorem above, as illustrated by Kargaev's counterexample, see \cite[vol. 1, p. 305]{KoosisLog}.

\ms\subsection{A hybrid theorem}
Beurling's and Levinson's Gap Theorems compliment each other by treating measures with sparse supports and fast decay
correspondingly. In this section we suggest a hybrid theorem that combines the features of both statements. In comparison
with Beurling's result it shows that the measure does not have to be zero on a long sequence of intervals, it just has to
be small on it. In regard to Levinson's theorem, our statement says that the measure does not have to decay fast along the whole
axis, just along a large enough set. One can show that  the statement is sharp in both scales.

\begin{thm}\label{BerLev}
Let $\mu$ be a finite measure on $\R$ whose Fourier transform vanishes on an interval. Suppose that there exists a sequence of disjoint intervals
$\{I_n\}$ such that

\begin{equation} \sum \frac{|I_n|\min\left(|I_n|,\log\frac 1{|\mu|(I_n)}\right)}{1+\dist^2(I_n,0)}=\infty.
\label{bleq}
\end{equation}
Then $\mu\equiv 0$.

\end{thm}

\begin{proof} We can assume that  $|I_n|\to\infty$ because any subsequence of intervals with uniformly bounded
lengths can be deleted from $\{I_n\}$ without affecting \eqref{bleq}.
Define the $|\mu|$-weight $W$ as
$$W=\left[|\mu|(I_n)(1+\dist^2(I_n,0))\right]^{-1}$$ on each $I_n$. If $\hat\mu$ vanishes on an interval then  for some $d>0$ there
exists a  $d$-uniform sequence $\L\subset\supp\mu$ satisfying \eqref{ur4}. Let
$$N=\{\ n\ |\ \#(\L \cap I_n)>\frac d2|I_n|\}.$$
 Note that the sequence $\{I_n\}_{n\not\in N}$ cannot be long
because otherwise $\L$ will not satisfy the density condition \eqref{density} on any short partition.
Therefore the part of the sum in \eqref{bleq} corresponding to $n\not\in N$ is finite and
$$\sum_{n\in \Z}\frac{\log W(\lan)}{1+\lan^2}\geqslant\sum_{\lan\in\cup_{k\in N} I_k}\frac{\log W(\lan)}{1+\lan^2}\gtrsim\sum_{n\in N}\frac{|I_n|\log \frac 1{|\mu|(I_n)}}{1+\lan^2}=\infty.$$

\end{proof}

\ms\subsection{De Branges' Gap Theorem}

\begin{thm}[de Branges, theorem 63 \cite{dBr}]\label{Branges}
Let $K(x)$ be a continuous function on $\R$ such that $K(x)\geqslant 1$, $\log K$ is uniformly continuous
and Poisson-unsummable.
Then there is no nonzero finite measure $\mu$ on $\R$ such that
\begin{equation}\int^{\infty}_{-\infty}Kd|\mu|<\infty
\label{messum}\end{equation}
and $\hat\mu$ vanishes on an interval.

\end{thm}

\begin{proof}
Suppose that $\mu$ satisfies \eqref{messum} and its Fourier spectrum has a gap.
Since $K$ is a $\mu$-weight there must exist a $d>0$ and a $d$-uniform sequence $\L\subset\supp \mu$
satisfying \eqref{ur4} with $K$ in place of $W$. Since $\L$ has
positive interior density and $\log K$ is uniformly continuous,
\eqref{ur4} implies that $\log K$ is Poisson-summable.
\end{proof}

\ms
\subsection{A theorem by Krein, Levinson and McKean}\label{secKLM}
Our next statement combines results by Krein  (part I in the statement below, case $p=2$) and
by Levinson and McKean  (part II, $p=2$).

\begin{thm}[Krein \cite{Krein1}, Levinson-McKean \cite{DM}]\label{KLM}
Let $\mu$  be a finite measure on $\R$, $\mu=w(x) dx$, where $w(x)\geqslant 0$.
Then

\ms\no I) If $\log w$ is Poisson-summable then for any $1\leqslant p\leqslant \infty$, $\GG^p_\mu=\infty$.

\ms\no II) If $\log w$ is monotone and Poisson-unsummable on a half-axis $(-\infty,x)$ or $(x,\infty)$ for some $x\in\R$ then for any $1< p\leqslant \infty$, $\GG^p_\mu=0$.

\end{thm}

\begin{proof}
If $\log w$ is Poisson-summable, denote by $H(z)$ the outer function in $\C_+$ satisfying $|H|=w$ on $\R$.
Then for any $a>0$ the measure $\eta=e^{-iax}\bar H(x)dx$ annihilates all exponentials with frequencies from $[0,a)$.
(Here we use the fact that the integral over $\R$ for any function from $H^1(\C_+)$ is 0.)
Since $|\eta|=\mu$, it follows that $\GG^p_\mu=\infty$ for any $1\leqslant p\leqslant \infty$.

\ms\no In the opposite direction, suppose that $\log w$ is  Poisson-unsummable and monotone on  $\R_+$. Consider a $\mu$-weight $W(x)=(w(x)(1+x^2))^{-1}$. If $\GG^p_\mu>2\pi d>0$, there exists a $d$-uniform sequence $\L$
satisfying \eqref{ur4}. Suppose that $\L$ satisfies \eqref{density} on a short partition $I_n=(b_n,b_{n+1}], \ b_0=0$.
Then, similarly to the  proof of theorem  \ref{Lev}, for some $0<c<1$, $cb_{n+1}<b_n$. Together with monotonicity
of $\log w=-\log W-\log(1+x^2)$, we obtain
$$\sum\frac{\log W(\lan)}{1+\lan^2}\gtrsim\sum_{n=0}^\infty\frac{\log W(b_n)|I_n|}{1+b_n^2}+\const$$$$
\gtrsim \sum_{n=0}^\infty\frac{\log W(cb_{n+1})|I_n|}{1+b_n^2}+\const\gtrsim \int_0^\infty \frac{-\log w(cx) dx}{1+x^2}+\const=\infty.$$

\ms\no
\end{proof}

\ms\subsection{A result by Borichev and Sodin on stability of type}\label{secBS}
If $I\subset \R$ is an interval and $D>0$ is a constant we denote
by $DI$ the interval concentric with $I$ of length $D|I|$.
Following \cite{BS}, for $\delta > 0$ and $x\in\R$, we denote
$$  I_{x,\delta} = [x - e^{-\delta |x|}, x +  e^{-\delta |x|}].$$
If $\mu$ and $\nu$ are two finite positive measures on $\R$ we write
$\mu\preccurlyeq\nu$ if there exist constants $\delta > 0,\ C > 0$, and $l > 0$, such that, for all $x\in R$,
$$\mu(I_{x,\delta} )\leqslant C(1 + |x|)^l \left(\nu(2I_{x,\delta}) + e^{-2\d |x|}\right)
.$$

\ms\no Instead of finite measures \cite{BS} deals with a wider class of polynomially growing measures
and uses the corresponding definition of type. As was mentioned in the introduction, in view of
the statements like lemma \ref{polynomial} above, such differences are not essential for the type problem
and the corresponding results are equivalent.

\begin{thm}\cite{BS}\label{BSthm}
If $\mu\preccurlyeq\nu$ then $\GG^2_\mu\leqslant\GG^2_\nu$.
\end{thm}

\begin{proof}
Let $\{a_n\}_{n\in\Z}$ be a strictly increasing discrete sequence of real points satisfying
$a_{-n}=-a_n$ and
$$a_{n+1}-a_n=2e^{-\d b_n},\ b_n=\frac{a_{n+1}+a_n}2\textrm{ for all }n\in\N,$$
where $\d>0$ is the constant from the definition of the relation $\mu\preccurlyeq\nu$.
Denote $I_n=(a_n,a_{n+1}]$. Let $W$ be a $\nu$-weight.
Then the step-function $W^*$ defined as
$$W^*(x)= 1+(1 + |b_n|)^{-l}\left[\frac 1{\nu(2I_n)+e^{-2\d b_n}}\int_{2I_n} Wd\nu\right]\textrm{ on each }I_n$$
is a $\mu$-weight, as follows from the condition $\mu\preccurlyeq\nu$. Assume that $\GG^2_\mu=2\pi d>0$.
Then there exists an $(d-\e)$-uniform sequence $\L\in\supp \mu$,
that satisfies \eqref{ur4} with $W^*$.
Our goal is to modify $\L$ into an $(d-\e)$-uniform sequence  in $\supp \nu$ satisfying \eqref{ur4} with $W$.

\ms\no Notice that WLOG we can assume that each interval
$2I_n$ contains at most one point of $\L$, see remark \ref{rem1}.
Choose $k_n$ so that $\lan\in 2I_{k_n}$.

\ms\no Now for each $\lan\in 2I_{k_n}$ choose a point $\alpha_n\in 2I_{k_n}\cap \supp \nu$ such that
$$W(\alpha_n)\leqslant \frac 1{\nu(2I_n)}\int_{2I_n} Wd\nu.$$
WLOG
$$\int_{2I_n} Wd\nu\geqslant e^{-\frac32\d b_n}$$
 for all $n$: otherwise we can increase the weight $W$ to satisfy
this condition and it will still remain $\nu$-summable. If $W$ is such a weight, then the interior density of the subsequence
of $\L$ that falls in the intervals $I_n$ satisfying $\nu(2I_n)\leqslant e^{-2\d b_n}$ must be zero: otherwise the sum \eqref{ur4} for $\L$ and $W^*$ would diverge. We can assume that $\L$ does not have
such points.
Then
$$ \log W(\alpha_n)\leqslant \log W^*(\lan)+2l\log(1+|\lan|)$$
and therefore
$A=\{\alpha_n\}$ satisfies \eqref{ur4} with $W$. By remark \ref{rem1}, $A$ has a $(d-\e)$-uniform subsequence. Hence $\GG^2_\nu\geqslant\GG^2_\mu-2\pi\e$.
\end{proof}

\ms\no Notice that our proof is $p$-independent, i.e. $\GG^2$ can be replaced with $\GG^p$ for any $1<p\leqslant\infty$ in the Borichev-Sodin result.

\ms\subsection{A sufficient condition by Duffin and Schaeffer}\label{secDS}
Our next statement is formulated in \cite{DS, BS} for Poisson-finite measures. Here we present
an equivalent finite version.

\begin{thm}\label{DSthm}
Let $\mu$ be  a finite positive measure on $\R$ such that
for any $x\in \R$
$$\mu([x-L,x+L])>c(1+x^2)^{-1}$$
for some $L,c>0$. Then $\GG^2_\mu\geqslant \pi/L$.

\end{thm}

\begin{proof}
If $\e>0$ consider $a_n=n(2L+\e)$. Then in every interval $(a_n-L,a_n+L)$ there exists
a subinterval $I_n$ of the length $\e$ satisfying
$$\mu(I_n)\geqslant \frac{d\e}{L(1+a_n^2)}.$$ It is left
to apply theorem \ref{SuffGen} to the sequence of centers of $I_n$.
\end{proof}

\ms\subsection{Benedicks' result on unions of intervals}\label{secBen}
 The following reslut contained in
\cite{Benedicks} provides  non-trivial examples of measures with positive
type. Until now, only a few examples of this kind existed in the literature.

\begin{thm}\label{ben}\cite{Benedicks}
Let $...<a_{-1}<a_0=0<a_1<a_2<...$ be a discrete sequence of points and let $I_n=(a_n,a_{n+1}]$
be the corresponding partition of $\R$.
Suppose that there exist positive constants $C_1,C_2,C_3$ such that

\ms\no 1) if  $$C_1^{-1}a_{2n+1}<a_{2k+1}<C_1a_{2n+1},$$
for some $n,k$, then
$$C_2^{-1}|I_{2n+1}|<|I_{2k+1}|<C_2|I_{2n+1}|;$$

\ms\no 2) for all $n$ $$C^{-1}_1|a_{2n+1}|<|a_{2n-1}|<C_1|a_{2n+1}| ;$$

\ms\no 3) for all $n$ $$|I_{2n+1}|>C_3\max (|I_{2n}|,1);$$

\ms\no 4)

\begin{equation}\sum\frac{|I_{2n+1}|^2}{1+a^2_{2n+1}}\left[\log_+\frac{|I_{2n+1}|}{|I_{2n}|}+1 \right]<\infty.\label{ur6}
\end{equation}

\ms
\no Then for any real number $A>0$
and $1\leqslant p<\infty$
there exists a nonzero function
$$f\in L^1(\R)\cap L^p(\R)\cap C^\infty(\R),\ \ \ \supp f\subset\cup I_{2n},$$
  such that $\hat f=0$ on $[0,A]$.

\end{thm}

\ms\no Here we will not concern ourselves with the condition $f\in C^\infty$. The rest of the statement,
i.e. the existence of $f\in L^1(\R)\cap L^p(\R)$, follows from theorem~\ref{main}. Moreover,
conditions 1 and 2 prove to be redundant.

\begin{proof}
Let $\{b_n\}_{n\in\Z}$ be a   sequence of positive integers, monotonically increasing to $\infty$
as $n\to\infty$ and as $n\to-\infty$, such that
if one replaces $|I_{2n+1}|$ in \eqref{ur6} with $b_{n}|I_{2n+1}|$ the series still converges.
Consider the sup-partition of $\{I_n\}$ defined in the following way.
Let
$$n_0=0, n_{k+1}-n_k=b_{n_k}$$
for $k> 0$ and
$$n_{k+1}-n_k=b_{n_{k+1}}$$
for $n<0$. Define $J_k=(a_{2n_k},a_{2n_{k+1}}]$. By 3,  the new partition satisfies the property
$|J_n|\to\infty$ and, because of monotonicity of $b_n$,
\begin{equation}\sum\frac{|J_{n}|^2}{1+\dist^2(0,J_n)}\left[\log_+\frac{|J_n|}{|J_n\cap\left(\cup I_{2k}\right)|}+1 \right]<\infty.\label{ur7}
\end{equation}
In particular $\{J_n\}$ is short.

\ms\no Let $C$ be a large positive number.
By $[\cdot]$ we will denote the integer part of a real number.
Define a sequence $\L$ as follows.
On each $J_k=(a_{2n_k},a_{2n_{k+1}}]$ place $N=[C|J_k|]$  points of $\L$ inside
$J_k\cap\left(\cup I_{2n}\right)$
so that
$$\lambda_{m_k}<\lambda_{m_{k+1}}<...<\lambda_{m_k+N}$$
and
$$|\left(\cup I_{2n}\right)\cap (a_{2n_k},\lambda_{m_k}]|= |\left(\cup I_{2n}\right)\cap (\lambda_{m_k+N},a_{2n_{k+1}}]|=|\left(\cup I_{2n}\right)\cap (\lambda_{l},\lambda_{{l+1}}]|,$$
for all $l, m_k\leqslant l<m_k+N-1$.

\ms\no Then conditions 3 and 4 of the theorem imply that $\L$ satisfies the energy condition \eqref{energy}
on $J_n$ and that $D_*(\L)=C$. Also
the measure
$$\nu=\charf_{\cup I_{2n}}\Pi$$
and $\L$ satisfy conditions of theorem \ref{SuffGen}. Therefore $\GG^p_\nu\geqslant 2\pi C$ for any $1\leqslant p\leqslant\infty$ which implies
the existence of the desired function $f$ satisfying $\hat f =0$ on $(0,2\pi C)$.

\end{proof}

\ms\no Notice that our proof actually produces $f\in L^\infty$.
If, in addition to the conditions of the theorem, $|I_{2n}|>\const>0$, then the remaining property
 $f\in C^\infty$ can be added with little effort. One would need to construct $f$
 supported on $\cup \frac 12 I_{2n}$ and then consider a convolution $f*\phi$ with a
 $C^\infty$-function $\phi$ with small support. In the general case  $f$ can
 be "smoothed out" using functions with  exponentially decreasing size of support and
 involving arguments like theorem \ref{BSthm}.

\ms
\section{Proofs: Auxiliary Statements}\label{lemmas}

\ms\no This section contains the results that will be needed to prove  theorems \ref{main} and
\ref{mainDiscrete}.

\ms\subsection{A measure with positive type}
The following  lemma is essentially proved, but not explicitly stated in \cite{GAP}.

\begin{lem}\label{disclemma}
Let $A=\{a_n\}$ be a discrete  sequence of distinct real numbers that has bounded gaps, i.e. $a_{n+1}-a_{n}<C$ for some $0<C<\infty$.
Denote by $b_n$ the middle of the interval $(a_n,a_{n+1})$,
$b_n=(a_n+a_{n+1})/2$.
Suppose that the sequence $A$
is $d$-uniform for some $d>0$.
Then there exists a finite positive measure supported on $B=\{b_n\}$,
$$\mu=\sum\beta_n\delta_{b_{n}},$$
satisfying
\begin{equation}0<\beta_{n}\leqslant \frac{\sqrt{a_{n+1}-a_{n}}}{1+a_n^2},\label{ur100}\end{equation}
such that $\GG^\infty_\mu\geqslant 2\pi d$.
\end{lem}

\begin{proof}
Let $\theta$ be the meromorphic inner function constructed for the sequence $A$ as in lemma 5 from \cite{GAP}.
By construction, the Clark measure $\nu=\mu_{-1}$ corresponding to $\theta$ is supported on $B$ and satisfies
\begin{equation}
 \nu(\{b_n\})\lesssim a_{n+1}-a_{n},\label{ur101}
\end{equation}
 see the estimate (7.3) in \cite{GAP}.

\ms\no Let $c=d-\e$.
As was proved in \cite{GAP}, if $\theta$ satisfies the conditions of lemma 5, \cite{GAP}, and $A$
is $d$-uniform,
 then
there exists $f\in K_\theta$ that is divisible by $e^{icz}$ in $\C_+$.
(This is one of the main steps in the proof of theorem 2, \cite{GAP}. See the part from the fourth line before claim 1
to the end of part I of the proof.)

\ms\no Then, by the Clark representation,
$2\pi i f=(1+\theta)Kf\nu$.
Since $1+\theta$ is outer, $Kf\nu$ is divisible by $e^{icz}$ in $C_+$.
Because $\e$ is arbitrary, by lemma \ref{polynomial}, the measure $\mu=|f|\nu/(1+x^2)$ satisfies $\GG^\infty_\mu\geqslant 2\pi d$.
Since $f\in L^2(\nu)$ and $\nu$ satisfies \eqref{ur101}, considering a constant multiple of  $\mu$  if necessary, we obtain \eqref{ur100}.
\end{proof}

\ms\subsection{Construction of an auxiliary sequence}
To apply our previous lemma in the main proofs we will need the following

\begin{lem}\label{sequences}
Let $B=\{b_n\}$ be a $d$-uniform sequence satisfying \eqref{density} and
\eqref{energy} on a short partition $\{I_n\}$.
Let $w(n)$ be a positive bounded function on $\Z$ such that
\begin{equation}\sum\frac {\log w(n)}{1+n^2}>-\infty.\label{w}\end{equation}

\ms\no Then for any $\e>0$ there exists a discrete
sequence $A=\{a_n\}$ satisfying:

\ms\no 1) $a_{n+1}-a_{n}<1/\e$.

\ms\no 2) Define the sequence $C=\{c_k\}$  as $c_k=\frac{a_{k+1}+a_k}2$. Then
the sequence $B'=B\cap C$ satisfies
$$\#\left(B'\cap I_n\right)\geqslant (d-e)|I_n|$$
for large enough $|n|$.

\ms\no 3) If $b_n=c_k$, i.e. $b_n$ is the middle of $(a_k, a_{k+1})$, then $a_{k+1}-a_k\leqslant w(n)$.

\ms\no 4) $A$ is $2d$-uniform

\ms\no 5)  $D^*(C\setminus B)\leqslant d+\e$.

\end{lem}

\begin{proof} Denote
$$l_n=\min (b_{n+1}-b_n, b_n-b_{n-1}, w(n)).$$
 Consider the sequence $P=\{p_n\}$  defined as
 $$p_{2n}=b_n-\frac 13 l_n, \ \ p_{2n+1}=b_n+\frac 13 l_n.$$
Choose a large $L>>1/\e$. Define the sequence
 $Q$ as follows:  if $p_{2n+2}-p_{2n+1}>L$, insert $M=[(p_{2n+2}-p_{2n+1})/L]$ points
 of $Q$ into the interval $(p_{2n+1},p_{2n+2})$ uniformly, i.e. at the points
 $$p_{2n+1}+k\frac{p_{2n+1}-p_{2n+2}}{M+1},\ k=1,2,...M.$$

 \ms\no Now put $A=P\cup Q$.


 \ms\no By our construction the  sequence $A$
 satisfies
 $$2\#(B\cap I_n)-2\leqslant \#(A\cap I_n)\leqslant 2\#(B\cap I_n)+\e|I_n|.$$
 To make $A$ satisfy the more precise density condition \eqref{density} with $2d$ we may
 need to delete some points of $B$ on each interval $I_n$ and consider a smaller
 sequence $B'$ in place of $B$ in the above construction. Note that we would have to delete
 at most $\e|I_n|$ points from $B$ on each $I_n$ and that $B'$ will satisfy the energy
 condition \eqref{energy} as a subsequence of $B$. After such an adjustment,
 $A$ will satisfy 1), 2), 3) and the density condition \eqref{density} with $2d$.

 \ms\no Note that $A$ satisfies the energy condition on  $\{I_n\}$. Indeed, let us denote $\Delta_n=\#(P\cap I_n)$ and $\Gamma_n=\#(Q\cap I_n)$.
  Then
 $$ \#(A\cap I_n)^2\log |I_n|-\sum_{a_n,a_k\in A\cap I_n}\log |a_n-a_k|=
 $$$$\left(\Delta_n^2\log |I_n|-\sum_{a_n,a_k\in P\cap I_n}\log |a_n-a_k|\right)+
 $$$$\left(\Gamma_n^2\log |I_n|-\sum_{a_n,a_k\in Q\cap I_n}\log |a_n-a_k|\right)+
 $$$$2\left(\Delta_n\Gamma_n\log |I_n|-\sum_{a_n\in P\cap I_n,a_k\in Q\cap I_n}\log |a_n-a_k|\right)=
 $$$$
 I+II+III.$$
  To estimate $I$ notice that for any $p_{2k}\in P\cap I_n$,
  $$  -\log (p_{2k+1}-p_{2k})\leqslant -\log w(k),$$
  by our choice of points $p_{2k},p_{2k+1}$.
  The rest of the terms in $I$ can be estimated by the similar terms for $B'$, i.e.
  $$I\lesssim \left(\#(B'\cap I_n)^2\log |I_n|-\sum_{b_n,b_k\in B'\cap I_n}\log |b_n-b_k|\right)$$
  $$-\sum_{p_{2k}\in P\cap I_n}\log w(k)+O(|I_n|^2).$$
  Since $B'$ satisfies the energy condition and because of \eqref{w} and shortness of the partition, $I$ will give finite
  contribution to the energy sum in \eqref{energy}.

 \ms\no To estimate $II$ notice that points in $Q$ are at a distance at least $L/2$ from each other. Therefore
 $$II\lesssim \left(\Gamma_n^2\log |I_n|-\sum_{0\leqslant n,k\leqslant \Gamma_n}\log |n-k|\right) + O(\Gamma_n^2)=
 $$$$\Gamma_n^2\log \frac{|I_n|}{\Gamma_n} + O(\Gamma_n^2)$$
 after estimating the sum via Stirling's formula. Notice that since $\Gamma_n<|I_n|$ and $$\log \frac{|I_n|}{\Gamma_n} <\frac{|I_n|}{\Gamma_n},$$
 the last quantity will also give finite contribution to \eqref{energy}.

 \ms\no Finally, $III$ can be estimated similarly to $II$. Just notice that any point $a_j$ in $P$ is at a distance at least $L/2$
 from $Q$ and therefore
 $$\Gamma_n\log |I_n|-\sum_{a_k\in Q\cap I_n}\log |a_j-a_k|\lesssim \Gamma_n \log \frac{|I_n|}{\Gamma_n}+O(|I_n|^2).$$
 Summing over all $a_j\in P\cap I_n$ and recalling that $\#(P\cap I_n)=\Delta_n\lesssim |I_n|$ we again get a finite quantity in \eqref{energy}.

\ms\no To prove 5), let us split $C$
into two subsequences:
$$C_1=\left\{(a_n+a_{n+1})/2\ |\  a_n, a_{n+1}\in P\right\}\textrm{ and }C_2=C\setminus C_1.$$
Notice that
$C_1\setminus B'$ has at most one point between each two points of $B'$. Therefore,
$$D^*(C_1\setminus B)\leqslant D^*(B)\leqslant d+\e.$$ Also,
if $2/L<<\e$ then $D^*(C_2)<\e$, because any two points of $C_2$ are at a distance at least $L/2$ from each other.
\end{proof}

\ms
\subsection{Existence of extremal measure with a spectral gap}
The lemma in this section can be viewed as a version of de Branges' theorem 66 from \cite{dBr}.
The last section of \cite{GAP} contains a discussion of that theorem and its equivalent reformulations.

\ms\no Here and throughout the rest of the paper we will use the standard notation $S(z)=e^{iz}$
for the exponential inner function in the upper half-plane. In general,
$S^a(z)=e^{iaz}$ is inner in $\C_+$ if $a>0$ and inner in $\C_-$ if $a<0$.

\begin{lem}\label{t66}
Let $\mu$ be a finite complex measure such that $\hat\mu\equiv 0 $ on $[0,a]$.
Let $W$ be a $|\mu|$-weight. Then there exists
a finite  measure
$\nu=\sum\alpha_n\delta_{\lan}$
concentrated on a discrete sequence $\L=\{\lan\}$ such that

\ms\no 1) $\L\subset\supp\mu$;

\ms\no 2)  $W$ is  a $|\nu|$-weight;

\ms\no 3) $\hat\nu\equiv 0$ on $[0,a]$;

\ms\no 4) The Cauchy integral $K\nu$   has no zeros in $\C$, $K\nu/S^a$ is outer in $\C_+$ and $K\nu$ is outer in $\C_-$.

\end{lem}

\begin{proof} It will be more convenient for us to assume that $\hat\mu\equiv 0$
on a symmetric interval $[-a,a]$. Then $\bar\mu$ has the same property. Hence
we can assume that the measure is real (otherwise consider $\mu\pm\hat\mu$).

\ms\no Consider the following set of measures on $\supp\mu$:
$$M_W=\{\ \nu\ |\ \int Wd|\nu|\leqslant 1,\ \hat\nu=0\textrm{ on }[-a,a], \ \supp\nu\subset\supp\mu,\ \nu=\bar\nu\}.$$
Notice that the set is non-empty, because $\mu\in M_W$, and convex. It is also $*$-weakly closed
in the space of all finite measures on $\supp\mu$.
Therefore by the Krein-Milman theorem it has an extreme point. Let $\nu$ be such a point.
We claim that it is the desired measure.

\ms\no First, let us note that $\hat\nu\equiv 0$
on  $[-a,a]$. It is well-known that this property is equivalent to the property that $\nu$ annihilates
the Payley-Wiener class $PW_a$, i.e. that for any bounded $f\in PW_a$,
$$\int fd\nu=0,$$
see for instance the last section of \cite{GAP}.

\ms\no Next, let us show that
the set of real
$L^\infty(|\nu|)$-functions $h$, such that $ \widehat{h\nu}\equiv 0$
on  $[-a,a]$, is one-dimensional and  therefore $h=c\in \R$.
(This is equivalent to the statement that the closure of $PW_a$ in $L^1(|\nu|)$ has deficiency 1, i.e. the space
of its annihilators is one dimensional.)

\ms\no Let there be a bounded real $h$ such that $ \widehat{h\nu}\equiv 0$
on  $[-a,a]$.
WLOG $h\geqslant 0$, since one can add constants, and $\int W|h|d|\nu|=1$. Choose
$0<\alpha<1$ so that $0\leqslant\alpha h<1$. Consider the measures
$\nu_1=h\nu$ and $\nu_2=(1-\alpha)^{-1}(\nu-\alpha\nu_1)$. Then both of them
belong to $M_W$ and $\nu=\alpha \nu_1+(1-\alpha)\nu_2$ which contradicts
the extremality of $\nu$.

\ms\no Now let us show that $\nu$ is discrete. Let $g$ be a continuous compactly supported real function on $\R$
such that $\int gd|\nu|=0$. By the previous part, there exists a sequence $f_n\in PW_a$,
$f_n\to g$ in $L^1(|\nu|)$. Indeed, otherwise there would exist a function $h\in L^\infty(|\nu|)$ annihilating all $f\in PW_a\cap L^1(|\nu|)$ and such that $\int hgd|\nu|=1$. Since $\int gd|\nu|=0$, $h\neq const$ and we would obtain a contradiction
with the property that the space of annihilators is one-dimensional.

\ms\no Since $\nu$ annihilates
$PW_a$ and $(f_n(z)-f_n(w))/(z-w)\in PW_a$ for every fixed $w\in\C\sm\R$,
$$0=\int \frac{f_n(z)-f_n(w)}{z-w}d\nu(z)=Kf_n\bar\nu(w)- f_n(w)K\nu(w)$$
and therefore
$$f_n(w)=\frac{Kf_n\nu}{K\nu}(w).$$
Taking the limit,
$$f=\lim f_n=\lim \frac{Kf_n\nu}{K\nu}=\frac{K g\nu}{K\nu}.$$

\ms\no Since all of $f_n$ are entire, one can show that the limit function $f$
is also entire. Indeed, first notice that there exists a positive function $V\in L^1(|\nu|)$
such that $f_{n_k}/V \to g/V$ in $L^\infty(|\nu|)$, for some subsequence $\{f_{n_k}\}$. To find such a $V$ first choose
$f_{n_k}$ so that $||f_{n_k}-g||_{L^1(|\nu|)}<3^{-k}$ and then put
$$V=1+\sum 2^k |f_{n_k}-g|.$$
Denote $F_k=f_{n_k}/V$ and $\eta=V|\nu|$. Then $F_k$ converge in $L^2(\eta)$ and by the Clark theorem
$(1-I)KF_k\eta$ converge in $H^2(\C_+)$,
where $I$ is the inner function whose Clark measure is $\eta$. Notice that
$$f_{n_k}=\frac{Kf_{n_k}\nu}{K\nu}=\frac{KF_k\eta}{K\nu}=\frac{(1-I)KF_k\eta}{(1-I)K\nu}.$$
Now let $T$ be a large circle in $\C$ such that $|(1-I)K\nu|>\const>0$ on $T$.
Denote $T_\pm=T\cap \C_\pm$ and let $m_T$ be the Lebesgue measure on $T$.
Since
$(1-I)KF_k\eta$ converge in $H^2(\C_+)$, $f_{n_k}$ converge in $L^1(T_+, m_T)$. Similarly,
$f_{n_k}$ converge in $L^1(T_-, m_T)$. By the Cauchy formula it follows that $f_{n_k}$ converge normally
inside $T$ and therefore $f$ is analytic inside $T$. Since such a circle $T$ can be chosen
to surround any bounded subset of $\C$, $f$ is entire.

\ms\no Since the numerator in the representation
$$f=\frac{K g\bar\nu}{K\bar\nu}$$
is analytic outside the compact support of $g$, the measure in the denominator
must be singular outside of that support: Cauchy integrals of non-singular measures have jumps at the real line on the support of the a.c. part, which would contradict the property that $f$ is entire. Choosing two different functions $g$ with disjoint supports
we conclude that $\nu$ is singular.

\ms\no Moreover, since $f$ is entire, the zero set of $f$ has to be discrete.
Since $\nu$ is singular, $K\nu$ tends to $\infty$ nontangentially in $\C_+$ at $\nu$-a.e. point and $f=0$ at $\nu$-a.e. point outside of the support of $g$.
Again, by choosing two different $g$ with disjoint supports, we can see that  $\nu$ is concentrated on a discrete set.

\ms\no It remains to verify 4).  Since we chose to deal with the symmetric
interval $[-a,a]$, we need to show that $K\nu/S^{\pm a}$ are outer
in $\C_\pm$ correspondingly.

\ms\no Let $J$ be the inner function corresponding to $|\nu|$ ($|\nu|$ is the Clark measure for $J$). Denote
$$G=\frac 1{2\pi i}(1-J)K\nu\in K_J.$$
As was mentioned in section \ref{Clark}, $G$ has non-tangential boundary values $|\nu|$-a.e. and
$$\nu=G|\nu|.$$
Since $K\nu$ is divisible by $S^a$ in $\C_+$, $G$ is divisible by $S^a$ in $\C_+$.
Suppose that  $G=S^a UH$ for some inner $U$.
Since the measure $\nu$ is real, $\bar G=G$, $|\nu|$-a.e.

\ms\no Let $F\in K_J$ be the function such that $\bar J G=\bar F$. Since $J=1$, $|\nu|$-a.e.,
$F=\bar G= G$, $|\nu|$-a.e. Since functions in $K_J$ are uniquely determined by their traces
on the support of the Clark measure $|\nu|$, $F=G=S^a UH$.
Notice that the function
$h=S^a (1+U)^2H$ also
belongs to $K_J$:
$$\bar J h=\bar J S^a (1+U)^2H=(\bar J G)\bar U(1+U)^2=\bar G\bar U(1+U)^2$$$$
=\overline{S^a (1+U)^2H}= \overline{h}\in \overline{H^2}(\C_+),$$
because $\bar U(1+U)^2$ is real a.e. on $\R$.
 Denote by $\gamma$ the measure from the Clark representation of $h$, i.e.
$$\gamma=h|\nu|, \ \ h=\frac 1{2\pi i}(1-J)K\gamma.$$
 Then
 $$\gamma=h|\nu|=\bar U(1+U)^2G|\nu|=\bar U(1+U)^2\nu.$$
  The Cauchy integral of $\gamma$ is divisible by $S^a$ in $\C_+$
 because $h$ is divisible by $S^a$ in $\C_+$. Since $\bar U(1+U)^2$ is real, a constant multiple of $\gamma$
 belongs to $M_W$.
 Since $U$ is non-constant and $|\nu|$ is the Clark measure for $J$, $\gamma$ is not a constant multiple
 of $\nu$.
 Again we obtain
a contradiction with the property that the space of annihilators is one-dimensional.

\ms\no Thus $G/S^a\in K_J$ is outer in $C_+$. Since $ J \bar G=\bar G$, the pseudocontinuation of $G$ does not have an inner factor
except $S^{-a}$ in $\C_-$ as well. Hence $K\nu/S^{\pm a}$ is outer in $\C_\pm$.

\ms\no If $G$ has a zero at $x=a\in \R$ outside of $\supp \nu$ then
$$\frac G{x-a}\in K_J$$
 and the measure
$$\gamma= \frac G{x-a}|\nu|$$
leads to a similar contradiction with the property that the space of annihilators is one-dimensional, since $(x-a)^{-1}$ is bounded and real on the support of $\nu$. Since
$G=\frac 1{2\pi i}(1-J)K\nu,$
 $K\nu$ does not have any  zeros on $\R$.
\end{proof}

\begin{rem}
A statement similar to lemma 9 from \cite{GAP}, where $S^a$ was replaced with an arbitrary inner function, can also be formulated
in the case of lemma \ref{t66}.
\end{rem}

\ms\subsection{Estimates of $\log|\theta|$ for a meromorphic inner function}

\begin{lem}\label{thetasum}
 Consider
a short partition $\{I_n\}$ of $\R$.
Consider the set of circles $T_n=\{z\ |\ |z-\xi_n|=2|I_n|\}$ where
$\xi_n\in I_n$.
Let $\theta$ be a meromorphic inner function such that
\begin{equation}\sum_n\frac{\#\left(\{\theta=1\}\cap 10I_n\right)|I_n|}{1+\dist^2(0,I_n)}<\infty.
\label{ur30}\end{equation}
Then the integrals
$$p_n=\int_{T_n}\left|\log |\theta(z)|\right|d|z|$$
satisfy
\begin{equation}\sum_n\frac {p_n}{1+\dist^2(0,I_n)}<\infty.\label{ur1}\end{equation}
\end{lem}

\begin{proof}
Suppose that $\theta=S^aB$
for some $a>0$ and some Blaschke product $B, \{B=0\}=\{a_n\}\subset \C_+$.
Then
$$\log |\theta|=\log |S^a| +\log |B|.$$
The integrals of $|\log |S^a||$ are summable because
$$|\log |S^a||\lesssim |I_n|$$
on $T_n$ and the sequence $\{I_n\}$ is short. To estimate the integral of $|\log |B||$
notice that
$$|\log |B(z)||=\sum_{a_k\in D_n}\left|\log\frac{|z-a_k|}{|z-\bar a_k|}\right|+\sum_{a_k\not\in D_n}\left|\log\frac{|z-a_k|}{|z-\bar a_k|}\right|,$$
where $D_n$ is the disk,  $D_n=\{z\ |\ |z-\xi_n|\leqslant 3|I_n|\}$.
Elementary estimates show that for any $a_n\in D_n$
$$\int_{T_n}\left|\log\frac{|z-a_k|}{|z-\bar a_k|}\right|d|z|\lesssim |I_n|.$$
Also, since for each $a_k\in D_n$ the argument of $\frac{z-a_k}{z-\bar a_k}$ increases by at least $\pi$ on the diameter of $D_n$, that is contained in  $10I_n$,
the number of points $a_k\in D_n$ is $\lesssim \#\left(\{\theta=1\}\cap 10I_n\right)$. Hence, because of \eqref{ur30},
such integrals will give a finite contribution to the sum in \eqref{ur1}.

\ms\no For each $a_k\not\in D_n$ one can show that
$$ \int_{T_n}\left|\log\frac{|z-a_k|}{|z-\bar a_k|}\right|d|z|\lesssim\int_{I_n} \frac {|I_n|\Im a_kdx}{(\Re a_k -x)^2+(\Im a_k)^2} .$$
Notice that
$$\sum_k\int_{I_n} \frac {\Im a_kdx}{(\Re a_k -x)^2+(\Im a_k)^2}=\int_{I_n} (\arg B)'\leqslant 2\pi\cdot\#\left(\{\theta=1\}\cap I_n\}\right)+\const.$$
Again, because of \eqref{ur30}, the integrals for $a_k\not\in D_n$ will give a finite contribution in \eqref{ur1}
\end{proof}

\ms\subsection{A version of the first BM theorem}
The following lemma is essentially a version of the so-called first Beurling-Malliavin theorem,
see also \cite{MIF2}.

\begin{lem}\label{BM1} Let $\{I_n\}$ be a long sequence of intervals and let $c$ be a positive constant.
Denote by $I_n'$ and $I_n''$ the intervals of the length $c|I_n|$ adjacent to  $I_n$ from the left and from the right correspondingly.
Let $u$ be a real function on $\R$ such
that
$$\Delta_n=\sup_{I_n''}u-\inf_{I_n'}u\geqslant d|I_n|$$
for all $n$ and for some $d>0$. Then $u$ is not a harmonic conjugate of a Poisson-summable function.
\end{lem}

\begin{proof}
Note that if $\ti u\in L^1_\Pi$ then $f=e^{-iu+\ti u}$ is an outer function in the Smirnov class in $\C_+$.
Moreover, $f$ belongs to the kernel $N^+[e^{iu}]$ of the Toeplitz operator with the symbol $e^{iu}$
in the Smirnov class. This contradicts a Toeplitz version of the first BM theorem, see section 4.4 of \cite{MIF2}.
\end{proof}

\ms\subsection{An estimate for an extremal discrete measure of positive type}
In this section we show that that a discrete measure of positive type, like in the statement
of lemma \ref{t66}, must have log-summable pointmasses. We start with the following elementary statement that can
be easily verified.

\begin{lem}\label{short2I} Let $\{I_n\}$ be a short sequence of intervals and let $C>1$. Denote
$$l_n=\sum_{I_m\cap CI_k\neq \emptyset}|I_m|.$$
 Then
$$\sum\frac{l_n|I_n|}{1+\dist^2(0,I_n)}<\infty.$$
\end{lem}

\ms\no Our main statement in this section is

\ms\begin{lem}\label{denlog}
Let $\nu$ be a  finite measure
$$\nu=\sum\alpha_n\delta_{\lan}$$
on a discrete sequence $\L=\{\lan\}$, such that
 $\alpha_n\neq 0$,
$\hat\nu\equiv 0$ on $[0,2\pi d]$, $K\nu$ does not have any zeros in $\C$, $K\nu/S^d$ is outer in $\C_+$ and  $K\nu$ is outer in $\C_-$.
 Then for any $\e>0$, $\L$ contains a $(d-\e)$-uniform subsequence and
\begin{equation}\sum\frac{\log |\alpha_n|}{1+n^2}>-\infty.\label{ur2}\end{equation}
\end{lem}

\begin{proof}
The statement that $\L$ contains a $(d-\e)$-uniform subsequence follows from the property that $\GG_\L\geqslant d$ and theorem \ref{mainGAP}.

\ms\no  To establish \eqref{ur2},
let us first show that there exists a short partition $\{I_k\}$  of $\R $ such that $\L$ satisfies \eqref{density} with
$d$
on that partition.

\ms\no Let $J$ be the inner function whose Clark measure is $|\nu|$. Then by the Clark theorem the function
$$Q(z)=(1-J)K\nu$$
belongs to $K_J$. It follows from the properties of $K\nu$ that $Q=S^dO$ in $\C_+$ for some outer $O$ and
$\bar JQ=\bar O$. Therefore the argument of $O$ satisfies $u=2\arg O=\arg J- dx$.
Notice that $\arg J$ is a growing function that is equal, up to a bounded term, to the counting function of
$\L$. Also, since $O\in H^2$, $\tilde u\in L^1_\Pi$.
If the desired short partition
$\{I_k\}$, where $\L$ satisfies \eqref{density}, does not exist then there exists a long sequence of intervals $\{J_k\}$ such that
\begin{equation}\left|\#(\L\cup J_k) - d|J_k|\right|\geqslant c_1 |J_k|\label{ur20}\end{equation}
for each $k$ and for some $c_1>0$.
First, let us assume that the difference in the left-hand side is positive for a long subsequence of $\{J_k\}$.
Let $J_k', J_k''$ denote the intervals of the length $c_2|J_k|, 0< c_2<<c_1$, adjacent to $J_k$ from the left
and from the right correspondingly.
Since $u'$ is bounded from below we get that
$$\Delta_k=\inf_{J_k''}u-\sup_{J_k'}u>c_3 |J_k|$$
for some $c_3>0$ on a long subsequence of $\{J_k\}$, if $c_2$ is small enough. By lemma \ref{BM1}, this contradicts the
property that $\tilde u\in L^1_\Pi$. If the difference in \eqref{ur20} is negative for a long subsequence
of $\{J_k\}$ then
lemma \ref{BM1} can be applied to $-u$ and the intervals $J_k', J_k''$ chosen so that
$J_k', J_k''\subset J_k,\ |J_k'|=|J_k''|=c_2|J_k|$, $J_k'$ shares its left endpoint with $J_k$ and $J_k''$ shares its right endpoint
with $J_k$, to arrive at the same contradiction. Hence a short partition where $\L$ satisfies \eqref{density} with $d$ does exist.

\ms\no Let  $\{I_k\}$  be such a partition.
 Let $\lambda_{n_k}\in I_k$ be such that
$$\log_-\alpha_{n_k}=\max_{\lan\in I_k} \log_-\alpha_{n}.$$
Suppose that \eqref{ur2} is not satisfied.
Then
\begin{equation}\sum_k \#(\L\cap I_k)\frac{\log_-\alpha_{n_k}}{1+n_k^2}\asymp \sum_k |I_k|\frac{\log_-\alpha_{n_k}}{1+n_k^2}=\infty.
\label{ur32}
\end{equation}

\ms\no Consider $\mu=\sum_n \delta_{\lan}$ the counting measure of $\L$. Since $\L$ satisfies \eqref{density}, $\mu$ is Poisson-finite.
Let $\theta$ be the inner function such that $\mu$ is its Clark measure. Since $\nu$ is finite, it can be represented as
$\nu=f\mu$ with $f\in L^2(\mu)$. Hence, by Clark theory,
$$F=\frac 1{2\pi i}(1-\theta)K\nu\in K_\theta$$
with $F(\lambda_n)=f(\lan)=\alpha_n$.

\ms\no For each $k$  consider the disk
$$D_k=\{z \ |\ |z-\lambda_{n_k}|< 2|I_k|\}$$
and its boundary circle $T_k=\partial D_k$.
Notice that for each $k$, $F$ does not have zeros in $D_k$. It does have poles at the points $\bar a_n\in \C_-$,
where $A=\{a_n\}$ are the zeros of $\theta$ in $\C_+$.
Hence in $D_k$ the function $F$ admits factorization $F= H_k/B_k$, where $B_k$ is the finite Blaschke product in $D_k$ ($|B_k|=1$ on $T_k$)
with zeros at $\bar A\cap D_k$, and $H_k$ is analytic without zeros in $D_k$.

\ms\no Notice that
\begin{equation}-\int_{T_k}\log_- |F(z)|d|z|\leqslant\int_{T_k}\log |F(z)|d|z|=\int_{T_k}\log |H_k(z)|d|z|\lesssim |I_k|\log\alpha_{n_k}
\label{ur2a}\end{equation}
by Jensen's inequality, because $F$ has only poles and no zeros in $D_k$. At the same time,
since $F\in K_\theta$, it belongs to $H^2(\C_+)$ and is equal to $\theta \bar G, G\in H^2(\C_+)$ in $\C_-$.
Denote by $T_k^\pm$ the upper and lower halves of $T_k$. Since the absolute value of an $H^2$ function
is bounded by
$$\const+\const\ |y|^{-1/2}$$
 inside the half-plane, we have
$$\int_{T_k}\log_+ |F(z)|d|z|\leqslant\int_{T^+_k}\log_+ |F(z)|d|z|+
\int_{T_k^-}\log_+| G(z)|d|z|
$$\begin{equation} +\int_{T_k^-}\log_+|\theta (z)|d|z|\lesssim |I_k| + v_k,\label{ur3}\end{equation}
where $\sum v_k/(1+a_{n_k}^2)<\infty$ by  lemma \ref{thetasum}, because
$$\#(\{\theta=1\}\cap 10I_k)\lesssim \sum_{I_m\cap 10I_k\neq\emptyset}|I_m|$$
and \eqref{ur30} is satisfied by lemma \ref{short2I}.

\ms\no Since $H_k\neq 0$ in $D_k$, $\log |H|$ is harmonic in $D_k$. Hence its values on
$I_k$ can be recovered from the values of $\log |H_k|=\log |F|$ on $T_k$ via the
Poisson formula. By \eqref{ur3}, the Poisson integral of $\log_+|F|$ will deliver
a small contribution, i.e. on each $I_k$ it will be equal to a function $h^+_k$ such that
$$\sum \int_{I_k}h^+_k(x)d\Pi(x)<\infty.$$
On the other hand, the Poisson integral of $\log_-|F|$ in $D_k$,
restricted to $I_k$,
will be equal to $h^-_k$, where $h^-_k(x)\asymp \log\alpha_{n_k}$ for all
$x\in I_k$ by \eqref{ur2a}.
Hence by \eqref{ur32}
$$\sum \frac{\int_{I_k}\log|H_k|dx}{1+\dist^2(0,I_k)}=-\infty.$$
Furthermore, similarly to the proof of lemma \ref{thetasum},
$$\deg B_k\lesssim \#(\{\theta=1\}\cap 5I_k)\lesssim \sum_{I_m\cap 5I_k\neq\emptyset} |I_m|.$$
Therefore by lemma \ref{short2I}
$$\sum_k \frac{|I_k|\deg B_k}{1+\dist^2(0,I_k)}<\infty.$$
Thus
$$\int_\R\log |F(x)|d\Pi\asymp \sum_k \frac{\int_{I_k}\left(\log |B_k(x)|+\log |H_k(x)|\right)dx}{1+\dist^2(0,I_k)}\lesssim
$$$$
\sum_k \frac{|I_k|\deg B_k+\int_{I_k}\log |H_k(x)|dx}{1+\dist^2(0,I_k)}=-\infty$$
and we obtain a contradiction.
\end{proof}

\begin{rem}
Using results of \cite{GAP} one can prove a slightly stronger statement that $\L$ itself is  $d$-uniform.
\end{rem}

\ms
\subsection{Equivalence of completeness in $L^p$ and $C_W$}
The  theorem we discuss in this section relates the type problem to Bernstein's study of
weighted uniform approximation,
see \cite{KoosisLog} or \cite{BS}.

\ms\no Consider a weight $W$, i.e. a lower semicontinuous function $W:\R\to [1,\infty)$
that tends to $\infty$ as $x\to\pm\infty$.
We define $C_W$ to be the space of all continuous functions on $\R$ satisfying
$$\lim_{x\to\pm\infty} \frac{f(x)}{W(x)}=0.$$
We define the norm in $C_W$ as
$$||f||=|| fW^{-1}||_\infty.$$

\ms\no The following is a well-known result by A. Bakan. For reader's convenience we supply a short proof.

\begin{thm}[\cite{Bakan}]\label{Bakan} Let $\mu$ be a finite positive measure on $\R$. Then
the system of exponentials $\EE_d$ is complete in $L^p(\mu)$ for some $1\leqslant p\leqslant \infty$ if and only if
there exists a $\mu$-weight $W\in L^p(\mu)$ such that
$\EE_d$ is complete in $C_W$.

\end{thm}

\begin{proof}
If $\EE_d$ is complete in $C_W$ for some $\mu$-weight $W\in L^p(\mu)$
then for any bounded continuous function $f$ there
exists a sequence $\{S_n\}$ of finite linear combinations of exponentials from $\EE_d$
such that $S_n/W$ converges to $f/W$ uniformly. Then $S_n$ converges to $f$ in $L^p(\mu)$.
Hence $\EE_d$ is complete in $L^p(\mu)$.

\ms\no Suppose that $\EE_d$ is complete in $L^p(\mu)$. Let $\{f_n\}_{n\in\N}$ be a  set
 of bounded continuous functions on $\R$, that is dense
 in $C_W$. Let $\{S_{n,k}\}_{n,k\in\N}$
be a family of finite linear combinations of exponentials from $\EE_d$ such
that
$$||f_n-S_{n, k}||_{L^p(\mu)}<4^{-(n+k)}.$$
Denote
$$W=1+\sum_{n,k\in \N}2^{n+k}|f_n-S_{n, k}|.$$
Notice that then $W\in L^p(\mu)$ and $S_{n, k}/W\to f_n/W$ uniformly as $k\to\infty$.
Since $\{f_n\}$ is dense in $C_W$, $\EE_d$ is complete in $C_W$.
\end{proof}

\ms
\section{Proofs of main results}\label{mainproofs}


\ms\subsection{Proof of theorem \ref{mainDiscrete}}$ $

\ms\no To prove that $\GG^p_\mu\geqslant 2\pi D$, WLOG we can assume that $B$ itself is a $d$-uniform sequence for some $d>0$
and that $w$ satisfies \eqref{ur10}.

\ms\no Fix a small $\e>0$. Let $C=\{c_n\}$ be the sequence provided by lemma \ref{sequences}. Then by lemma \ref{disclemma} (applied to
$C$ and $w^2$)
there exists a finite positive measure $\nu=\sum \sigma_n\delta_{c_n}$ concentrated on $C$,
satisfying
$$0<\sigma_n<w(k)\textrm{ for  } c_n=b_k\textrm{ and }\GG^\infty_\nu\geqslant 2\pi (2d).$$

\ms\no Let $\theta$ be the Clark inner function corresponding to $\nu$. Then there exists  a function  in $K_\theta$ divisible
by $S^{2\pi (2d-\e)}$ in the upper half-plane, i.e. $S^{2\pi (2d-\e)}h\in K_\theta$ for some $h\in H^2$:
if $\widehat{\phi\nu}=0$ on $[0,2\pi (2d-\e)]$ for some $\phi\in L^\infty(\nu)$, put
$$h=\frac 1{2\pi i}(1-\theta)K\phi\nu.$$

 \ms\no By lemma \ref{sequences}, $D^*(C\setminus B)< d+\e$.
Let $J$ be an inner function such that $\{J=1\}=C\setminus B$.
By a version
of the Beurling-Malliavin theorem, see \cite{MIF1} section 4.6, the kernel of the Toeplitz operator with the symbol
$S^{2\pi (-d-\e)}J$ in $H^\infty$ is non-empty, i.e. there exists a function  $g\in H^\infty(\C_+)$ such that
$$S^{2\pi (-d-\e)}Jg\in \bar H^\infty.$$
Since
$$\bar\theta S^{2\pi (-d-\e)}Jg S^{2\pi (2d-\e)}h=\bar\theta S^{2\pi (d-2\e)}Jgh\in \bar H^2,$$
we have
$$S^{2\pi (d-2\e)}Jgh\in K_\theta.$$
 Since $K_\theta$ is closed under division by inner components,
$S^{2\pi (d-2\e)}gh\in K_\theta$ and therefore
$$p=S^{2\pi (d-2\e)}Jgh-S^{2\pi (d-2\e)}gh=S^{2\pi (d-2\e)}(J-1)gh\in K_\theta.$$

\ms\no By the Clark representation formula, $p=\frac1{2\pi i}(1-\theta)Kp\nu$, and since $1-\theta$ is outer, $Kp\nu$ is divisible by $S^{2\pi (d-2\e)}$
in $\C_+$. Notice that $p=(1-J)gh=0$ on $C\setminus B=\{J=1\}$ and $p\in L^\infty(\nu)$ on $B\cap C$. Therefore, if $\eta$ is the restriction of $\nu$
on $B\cap C$, the existence of such $p$ implies
$$\GG^\infty_\eta\geqslant 2\pi (d-2\e).$$

\ms\no For any $\e>0$, the measure $\eta$ constructed as above will have a bounded density with respect to
$\mu$. Hence
$\GG^\infty_\mu\geqslant 2\pi d.$

\ms\no To prove the second part of the statement suppose that $\log(|n_B|+1)\in L^1_\Pi$ but $\GG^p_\mu>2\pi A> 2\pi D$ for
some $A>D$. WLOG assume that the counting function $n_B$ is non-zero outside of $[-1,1]$ and define
$$W(b_k)=\begin{cases}\frac1{2^{n}n_B(2^n)w(k)}+1\textrm{ if }b_k\in (2^{n-1},2^n]\textrm{ for some }n\in\N\\
-\frac1{2^{n}n_B(-2^n)w(k)}+1\textrm{ if }b_k\in (-2^{n},-2^{n-1}]\textrm{ for some }n\in\N\\
1\textrm{ if }b_k\in (-1,1]\end{cases}.
$$
 Then $W$ is a $\mu$-weight and by lemmas \ref{t66}  and \ref{denlog} there exists a measure $\nu=\sum \alpha_k\delta_{b_{n_k}}$ supported on $B'=\{b_{n_k}\}\subset B$ such that
 $W$ is a $\nu$-weight, $B'$ is an $A$-uniform sequence and $\alpha_k$ satisfy \eqref{ur2}. Since $W$ is a $\nu$-weight, $|\alpha_k|\leq C/W(b_{n_k})$. Since
 $\alpha_k$ satisfy \eqref{ur2}, the definition of $W$ implies that
 $$\sum_k\frac {\log w(n_k)}{1+{n_k}^2}>-\infty.$$
 Hence $D\geqslant A$ and we obtain a contradiction.
$\square$

\ms\subsection{Proof of theorem \ref{main}}
$ $

\ms\no I) First, suppose that $\GG^p_\mu\geqslant a$ for some $1<p\leqslant \infty$. Then for any $d>0,\ 2\pi d<a$, there exists $f\in L^p(\mu)$ such that
$\widehat{f\mu}=0$ on $[0,2\pi d]$. Let $W$ be a $\mu$-weight. Denote $V=W^{1/q}$ where $\frac1p+\frac1q=1$. Then
$$\int V|f|d\mu<\infty.$$
Therefore by lemma \ref{t66} there exists a discrete measure $\nu=\sum \alpha_n\delta_{\lan},\ \L=\{\lan\}\subset\supp\mu$
such that $\hat\nu=0$ on $[0,2\pi d]$, $V$ is a $|\nu|$-weight and $\nu$ satisfies the rest of the conditions of lemma
\ref{denlog}. Then by lemma \ref{denlog}, $\L$ contains a  $d$-uniform subsequence $\L'$ and $\alpha_n$ satisfy \eqref{ur2}.
Since
$V$ is a $\nu$-weight, $V(\lan)<C/|\alpha_n|$ for all $n$. It is left to notice
that $\log W(\lan)=q\log V(\lan)$ and therefore
$$\sum_{\lan\in \L'}\frac{\log W(\lan)}{1+\lan^2}<\infty.$$
II) Now suppose that
$\GG^p_\mu< d<a$ for some $1<p\leqslant \infty$.

\ms\no Since $\GG^p_\mu< d$, by theorem \ref{Bakan} there exists a $\mu$-weight $W$ such that
finite linear combinations of exponentials from $\EE_{d-\e}$ are dense in $C_W$ for some $\e>0$.
Suppose that there exists a $d$-uniform sequence
$\L=\{\lan\}\subset \supp \mu,$  satisfying
\eqref{ur4}.
Then by theorem \ref{mainDiscrete} there exists a measure $\nu=\sum\alpha_n \delta_{\lan}$
such that $|\alpha_n|\leqslant W^{-1}(\lan)/(1+\lan^2)$ and $\hat\nu=0$ on $[0,d-\e]$. Then
the finite measure $W\nu$ annihilates all functions $e^{ict}/W,\ c\in [0,d-\e]$.  This contradicts
completeness of $\EE_{d-\e}$ in $C_W$.
$\square$

\bs\bs

\end{document}